\documentclass[lettersize,journal]{IEEEtran}
\IEEEoverridecommandlockouts

\usepackage{cite}
\usepackage{amsmath,amssymb,amsfonts}
\usepackage[capitalise]{cleveref}
\usepackage[]{footmisc}
\usepackage{xcolor}
\usepackage{algorithm}
\usepackage{algorithmic}
\usepackage[algo2e]{algorithm2e} 
\usepackage{graphicx}
\usepackage{textcomp}
\usepackage{psfrag}
\usepackage{epsfig}
\usepackage{subcaption}
\usepackage{float}
\usepackage{xcolor}
\usepackage{empheq}
\usepackage{mathtools}
\usepackage[numbers,sort&compress]{natbib}
\usepackage{tabularx}
\usepackage{booktabs}
\usepackage{diagbox}
\usepackage{multicol}
\newtheorem{mydef}{Definition}

\newtheorem{myrem}{Remark}
\newtheorem{mylem}{Lemma}
\newtheorem{mythr}{Theorem}

\newtheorem{myproblem}{Problem}

\begin{document}

\title{\LARGE Sparse $H_\infty$ Controller  for Networked Control Systems: Non-Structured and Optimal Structured Design}

\author{Zhaohua~Yang,~Pengyu~Wang,~Haishan~Zhang,~Shiyue~Jia,~Nachuan~Yang,\\~Yuxing~Zhong~and~Ling~Shi,~\IEEEmembership{Fellow,~IEEE}
	\thanks{\rm Z. Yang, P. Wang, H. Zhang, S. Jia, Y. Zhong and L. Shi are with the Department of Electronic and Computer Engineering, Hong Kong University of Science and Technology, Clear Water Bay, Hong Kong SAR. L. Shi is also with the Department of Chemical and Biological Engineering, Hong Kong University of Science and Technology. (email: zyangcr@connect.ust.hk; pwangat@connect.ust.hk; hzhangdo@connect.ust.hk; sjiaac@connect.ust.hk; yuxing.zhong@connect.ust.hk; eesling@ust.hk).}
	\thanks{N. Yang is with the Department of Electrical and Computer Engineering, University of Alberta,
Edmonton, Alberta T6G 1H9, Canada (email:nachuan1@ualberta.ca).}
    \thanks{The work is supported by the Hong Kong RGC General Research Fund 16203723 (Corresponding author: Yuxing Zhong).}
  }
\maketitle

\begin{abstract}
This paper provides a comprehensive analysis of the design of optimal structured and sparse $H_\infty$ controllers for continuous-time linear time-invariant (LTI) systems. 
Three problems are considered.
First, designing the sparsest $H_\infty$ controller, which minimizes the sparsity of the controller while satisfying the given performance requirements. Second, designing a sparsity-promoting $H_\infty$ controller, which balances system performance and controller sparsity. Third, designing a $H_\infty$ controller subject to a structural constraint, which enhances system performance with a specified sparsity pattern. For each problem, we adopt a linearization technique that transforms the original nonconvex problem into a convex semidefinite programming (SDP) problem. Subsequently, we design an iterative linear matrix inequality (ILMI) algorithm for each problem, which ensures guaranteed convergence. We further characterize the first-order optimality using the Karush-Kuhn-Tucker (KKT) conditions and prove that any limit point of the solution sequence generated by the ILMI algorithm is a stationary point. For the first and second problems, we validate that our algorithms can reduce the number of non-zero elements and thus the communication burden through several numerical simulations. For the third problem, we refine the solutions obtained in existing literature, demonstrating that our approaches achieve significant improvements.
\end{abstract}


\IEEEpeerreviewmaketitle
\section{Introduction}
\IEEEPARstart{N}{etworked} control systems (NCS) have found ubiquitous industrial applications and have received much attention from the research community. In such systems, information is exchanged via a communication network among system components such as controllers, sensors, and actuators \cite{sandberg2015cyberphysical}. Such types of systems have found applications in power networks \cite{teixeira2010networked}, transportation networks \cite{ran2012modeling}, and sensor and actuator networks \cite{verdone2010wireless}. Within NCS, we can design a centralized feedback control strategy that requires each controller to access the full state. Although this approach achieves optimal control performance, it typically imposes a prohibitively high communication burden in large-scale systems. Moreover, in practical applications, the communication and computation capabilities are often limited \cite{gupta2015survey}, which significantly restricts the number of communication channels and, consequently, the scalability and potential of NCS with such centralized designs. This motivates us to design control inputs using partial state information, which corresponds to a sparse controller gain matrix. 

Extensive work has been conducted on the optimal sparse controller design. In the seminal work \cite{lin2013design}, the authors proposed an alternating direction method of multiplier (ADMM) algorithm for designing sparse controllers. Instead of ADMM, in \cite{fardad2014design}, the authors proposed an iterative convex programming algorithm. More recently, various algorithms have been proposed for the structured and sparse controller design \cite{yang2024log,dhingra2016method,babazadeh2016sparsity,cho2024iterative,yang2025sparse}. The sparsity mentioned above refers to the elementwise sparsity over the gain matrix. In contrast, in the landmark work \cite{polyak2013lmi}, the authors proposed the concept of row sparse and column sparse, which motivates the sensor and actuator selection problems. In such types of problems, a change of variable technique can be applied to turn the feasible region into a convex set. In \cite{dhingra2014admm,zare2018optimal,zare2019proximal}, the authors proposed various efficient algorithms in large-scale problems, including ADMM, proximal gradient, and quasi-Newton. The sparse control has also been investigated in topology design \cite{ding2019sparsity}, consensus network \cite{lin2012identification}, covariance completion \cite{zare2015alternating}, target tracking \cite{masazade2012sparsity}, hands-off control \cite{nagahara2015maximum,nagahara2020sparsity}, observer design \cite{yang2024sparsity}, and remote state estimation \cite{zhong2024sparse}. However, to the best of our knowledge, most previous works focus on $H_2$ performance, with very few considering $H_\infty$ performance, which is another important performance criterion. Designing structured and sparse $H_\infty$ controllers remains an challenging problem due to its nonconvexity. Additionally, we cannot design a first-order
algorithm for analysis since we lack an analytical expression
for a given controller. In \cite{ferrante2020lmi}, the authors designed a structured $H_\infty$ controller by considering a convex subset of the nonconvex feasibility domain, which omitted some potential solutions and could lead to suboptimality. Therefore, more advanced algorithms capable of thoroughly exploring the nonconvex feasibility region remain to be developed, which motivates our research in this paper.

In this paper, we provide a systematic solution to the optimal structured and sparse $H_\infty$ controller design. The main contributions are multifold.
\begin{enumerate}
    \item We consider three typical scenarios:
    \begin{enumerate}
        \item sparse controller design with bounded $H_\infty$ performance: minimizing the controller's sparsity while satisfying a specified performance threshold.
        \item sparsity-promoting $H_\infty$ controller design: exploring the trade-offs between system performance and controller sparsity.
        \item structured $H_\infty$ controller design: enhancing the system performance subject to specified structural constraints on the controller. 
    \end{enumerate}
    For a) and b) we, for the first time, design sparse $H_\infty$ controllers without structural constraints. For c), we refine the solutions obtained in \cite{ferrante2020lmi}. Using the same numerical examples, our simulations demonstrate that our approaches achieve significant improvements over the solutions presented in \cite{ferrante2020lmi}.
    \item We develop a novel linearization technique that transforms the problems mentioned above into convex SDP problems, allowing for the explicit incorporation of structural constraints and $l_1$ regularization terms. We further show that any feasible solution to the transformed problem is also feasible for the original problem, motivating an ILMI algorithm to solve each problem.
    \item By incorporating proximal terms into the objective function, we demonstrate that the algorithm exhibits sufficient decrease and is guaranteed to converge. Moreover, we characterize the first-order optimality using KKT conditions, showing that any limit point of the solution sequence satisfies the KKT conditions of the original problem and is thus a stationary point.
\end{enumerate}

This paper is organized as follows. In \cref{Preliminaries}, we provide background on $H_\infty$
  performance, sparsity control, and existing structured $H_\infty$
  controller design methods, and formulate the problem. In \cref{Main Results}, we separately consider the three different problems, develop the linearization technique, present the optimality conditions, and propose the ILMI algorithm. In \cref{Simulations}, we provide numerical simulations demonstrating the effectiveness of our methods for scenarios (i) and (ii), as well as the superiority of our approaches in scenario (iii). Finally, in \cref{Conclusions}, we conclude the paper.

\textit{Notations}:
$$
\begin{aligned}
&\begin{array}{ll}
\mathbb{R}^{m\times n} & \text{real matrix with dimension} \ m\times n\\
\mathbb{S}^{n} & \text{symmetric matrix with dimension} \ n\times n\\
X^\top\ (X^H)& \text{transpose (conjugate transpose) of}\ X\\
X_{ij}\ ([X]_{ij}) & \text{element at $i$-th row and $j$-th column}\\
X_k\ ([X]_k) & \text{matrix sequence at $k$-th iteration} \\
X>0\ (X\ge0) & \text {element-wise positive (non-negative)}\\
X\succ 0\ (X\succeq 0) &  \text {positive definite (semidefinite)}\\
\mathrm{Tr}(X) & \text{trace of}\ X\\
I, 0, \mathbf{1} & \text{(identity, zero, all-ones) matrix}\\
\sigma_{\max} & \text{maximum singular value of a matrix} \\
\mathrm{Sym}(A) & A+A^\top \\
\mathrm{abs}(\cdot) & \text{element-wise absolute value operation}\\
\langle X, Y\rangle & \text{Frobenius inner product}\ (\mathrm{Tr}(X^\top Y))\\
\circ & \text{Hadamard product}\\
||\cdot||_F & \text{Frobenius norm}\\
\left[X_1|X_2|\ldots|X_n\right] & \text{row-stacked matrix}\\
\mathrm{Re}(\cdot) & \text{real part}\\
\triangleq & \text{equality for definition}\\
\otimes & \text{Kronecker product} \\
\left[\begin{smallmatrix}A & * \\ B^\top &C\end{smallmatrix}\right] & \left[\begin{smallmatrix}A & B \\ B^\top &C\end{smallmatrix}\right]
\end{array}
\end{aligned}
$$
\section{Preliminaries}\label{Preliminaries}
\subsection{$H_{\infty}$ Performance}
In this paper, we consider a continuous-time LTI system given by
\begin{equation}\label{original system}
\begin{split}
    &\dot{x}(t) = Ax(t) + Bu(t) + Gd(t),\\
    &z(t) = Cx(t) + Du(t) + Hd(t),
\end{split}
\end{equation}
where $x\in \mathbb{R}^{n_x}$ is the system state, $u\in \mathbb{R}^{n_u}$ the control input, $d\in \mathbb{R}^{n_d}$ the exogenous disturbance, and $z\in \mathbb{R}^{n_z}$ the controlled output. It is standard to assume $(A, B)$ is stabilizable. We adopt a static state feedback control law to regulate the system \eqref{original system}, i.e.,
\begin{equation}\label{state feedback control law}
    u(t) = Kx(t),
\end{equation}
where $K\in \mathbb{R}^{n_u\times n_x}$ is the gain matrix to be designed. We can rewrite system \eqref{original system} as follows
\begin{equation}\label{reformulated system}
\begin{split}
    &\dot{x}(t) = (A+BK)x(t) + Gd(t),\\
    &z(t) = (C+DK)x(t) + Hd(t).
\end{split}
\end{equation}
Therefore, the influence of the exogenous disturbance $d$ on the controlled output $z$ is described by 
\begin{equation}\label{relationship between z and d}
    z(s) = T_{zd}(s)d(s),
\end{equation}
where $T_{zd}(s)$ is the transfer function given by
\begin{equation}\label{transfer function}
    T_{zd}(s) = (C+DK)(sI-(A+BK))^{-1}G + H.
\end{equation}
The $H_{\infty}$ norm of the system \eqref{reformulated system} is represented as $||T_{zd}(s)||_\infty$, where $||T_{zd}(s)||_\infty:=\mathrm{max}_{\omega\in \mathbb{R}}\sigma_{\max}\{T_{zd}(j\omega)\}$. It is well-known that the following result holds \cite{dullerud2013course}
\begin{equation}
    \int_0^{\infty}z^H(t)z(t)dt \leq ||T_{zd}(s)||^2_\infty \int_0^{\infty}d^H(t)d(t)dt,
\end{equation}
which shows that $||T_{zd}(s)||^2_\infty$ represents the largest energy amplification from the exogenous disturbance to the system output. The following fundamental result provides necessary and sufficient conditions for $H_\infty$ controller synthesis with a given performance threshold.
\begin{mylem}[\cite{iwasaki1994all}]\label{h infinity lemma}
    The system \eqref{reformulated system} is asymptotically stable (i.e., for all eigenvalues $\lambda$ of $A+BK$, $\mathrm{Re}(\lambda)<0$) and $||T_{zd}(s)||_\infty \le \gamma$ if and only if there exists a symmetric positive definite matrix $P\succ0$ such that 
    \begin{equation}
    \label{lemma1}
    \left[ \begin{matrix}\centering
	&\mathrm{Sym}(P(A+BK))		&PG		&(C+DK)^\top		\\
	&G^\top P		&-\gamma I		&H^\top		\\
	&C+DK		&H		&-\gamma I		\\
    \end{matrix} \right] \preceq 0.
    \end{equation}
\end{mylem}
The optimal $H_\infty$ norm is the minimal $\gamma$ such that \eqref{lemma1} is feasible. Therefore, we can obtain the optimal $H_\infty$ norm via the following SDP problem.

\begin{subequations}\label{unconstrained h infinity}
\begin{align}
\mathop{\min}_{X,Y,\gamma} &\ \gamma\\
\mathrm{ s.t. } &\ X \succ 0, \\
& \left[ \begin{matrix}
	       &\mathrm{Sym}(AX+BY)		&G		&(CX+DY)^\top		\\
	       &G^\top		&-\gamma I		&H^\top		\\
	       &CX+DY		&H		&-\gamma I		\\
            \end{matrix} \right] \preceq 0.
\end{align}
\end{subequations}

The optimal $H_\infty$ controller can be recovered by $K^*=Y^*{X^*}^{-1}$ with the corresponding global minimal $||T_{zd}(s)||_\infty=\gamma^*$, where $(Y^*,X^*,\gamma^*)$ is the solution to Problem \eqref{unconstrained h infinity}.
\subsection{Sparsity in NCS}
The feedback control strategy typically constructs a communication network between the control input and the system state. This network is generally dense, requiring each local controller to access the states of all subsystems. However, this could impose a prohibitively high communication burden on large-scale systems, which often have limited computation and communication capabilities. Therefore, it is preferable to design the control input that utilizes local state information in large-scale systems. This limited information exchange is reflected in the sparsity of the gain matrix. 
\subsubsection{Sparse Controller and $l_1$ Relaxation}
Considering the state feedback control law, each local controller can be computed as 
\begin{equation}
    \forall  1\le i \le n_u, [u(t)]_i = \sum_{1\le j\le n_x} [K]_{ij}[x(t)]_j.
\end{equation}
If $[K]_{ij}\neq0$, $[u(t)]_i$ needs access to $[x(t)]_j$, and otherwise $[u(t)]_i$ does not. Therefore, it is natural to define the communication burden as the number of nonzero elements of the matrix $K$. This is exactly the $l_0$ norm, represented as $||K||_0$. However, the $l_0$ norm is non-convex and discontinuous, and thus cannot be tackled directly. A common and effective technique to address this difficulty is to alternatively consider the $l_1$ norm defined as 
\begin{equation}
    ||K||_1 = \sum_{1\le i\le n_u}\sum_{1\le j\le n_x}|K_{ij}|.
\end{equation}
This $l_1$ relaxation turns the non-convex $l_0$ norm into a convex one and paves the way for our later analysis.

\subsubsection{Structured Controller}
A structured controller refers to a controller subject to a structural constraint $K\in \mathcal{S}$, where $\mathcal{S}$ is a specified linear subspace. The corresponding structured identity matrix is denoted as $I_\mathcal{S}$ defined by
\begin{equation}
    [I_\mathcal{S}]_{ij} \triangleq \begin{cases}
            1, &\text{if} \ K_{ij}\ \text{is a free variable},\\
            0, &\text{if} \ K_{ij} = 0\ \text{is required}.
        \end{cases}
\end{equation}
Then we define the complementary structured identity matrix as $I_{\mathcal{S}^c} \triangleq \mathbf{1} - I_\mathcal{S}$. Therefore, the structural constraint $K\in \mathcal{S}$ can be alternatively expressed as $K\circ I_{\mathcal{S}^c} = 0$.


\subsection{A Convex Approach for Structured Control}
Designing a structured $H_\infty$ controller is typically NP-hard \cite{blondel1997np}, meaning we cannot equivalently transform this problem into solving an SDP feasibility problem. To tackle this difficulty, Ferrante et al. \cite{ferrante2020lmi} proposed exploring solutions within a convex subset of the nonconvex feasibility region. Relevant results are provided below.
\begin{mythr}[\cite{pipeleers2009extended}]\label{equivalence}
    Given $P\succ0$ and $K\in \mathbb{R}^{n_u\times n_x}$, the following items are equivalent.
 \begin{equation}
        (i)\begin{bmatrix}
            &\mathrm{Sym}((A+BK)P) &G &P(C+DK)^\top\\
            &G^\top &-\gamma I &H^\top\\
            &(C+DK)P &H &-\gamma I
        \end{bmatrix}\prec0;\end{equation}\\
        ($ii$) There exists $X\in \mathbb{R}^{n_x\times n_x}$ and $\alpha>0$ such that:
        \begin{equation}
        \scalebox{0.87}{$\mathrm{Sym}\left(\begin{bmatrix}
                &(A+BK)X &\alpha(A+BK)X+P &0 &G\\
                &-X  &-\alpha X  &0  &0\\
                &(C+DK)X  &\alpha(C+DK)X  &-\frac{\gamma}{2}I  &H\\
                &0  &0  &0  &-\frac{\gamma}{2}I
            \end{bmatrix}\right)\prec 0$}.
        \end{equation}

\end{mythr}
\begin{mythr}[Theorem 2 of \cite{ferrante2020lmi}]\label{ferrante results}
    Let $\{S_1,S_2,\ldots,S_k\}$ be a basis of $\mathcal{S}$ and $$L\triangleq [S_1|S_2|\ldots|S_k].$$ Define the following \textit{structured set} $$\Upsilon  \triangleq\{Q\in \mathbb{R}^{n_x\times n_x}: \exists \Lambda\in \mathbb{S}^k\ \mathrm{s.t.} \ L(I_k\otimes Q)=L(\Lambda\otimes I_{n_x})\}.$$
    Given $\gamma>0$, if there exists $P\succ 0$, $R\in \mathrm{span}\{S_1,S_2,\ldots,S_k\}$, $\alpha>0$, and $X\in \mathbb{R}^{n_x\times n_x}$ such that:
    \begin{subequations}\label{ferrante feasibility}
        \begin{align}
            &X\in \Upsilon\label{ferrante con1},\\
            &\scalebox{0.9}{$\mathrm{Sym}\left(\begin{bmatrix}
                &AX+BR &\alpha(AX+BR)+P &0 &G\\
                &-X  &-\alpha X  &0  &0\\
                &CX+DR  &\alpha(CX+DR)  &-\frac{\gamma}{2}I  &H\\
                &0  &0  &0  &-\frac{\gamma}{2}I
            \end{bmatrix}\right)\prec 0$}.\label{ferrante con2}
        \end{align}
    \end{subequations}
    Then, $X$ is nonsingular and $K=RX^{-1}\in \mathcal{S}$ and $||T_{zd}(s)||_\infty \le \gamma$.
\end{mythr}

The corresponding problem for a minimal value of $\gamma$ can be expressed as follows \cite{ferrante2020lmi}:
\begin{subequations}\label{ferrante minimal gamma}
\begin{align}
\mathop{\min}_{P,X,R,\alpha,\gamma} &\ \gamma\\
\mathrm{ s.t. } &\ \eqref{ferrante con1},\eqref{ferrante con2}, P\succ 0, \alpha>0,\\
&\ R\in \mathrm{span}\{S_1,S_2,\ldots,S_k\}.
\end{align}
\end{subequations}
A structured $H_\infty$ controller can be computed as $\bar{K}=\bar{R}\bar{X}^{-1}\in \mathcal{S}$ with the corresponding $H_\infty$ norm $||T_{zd}(s)||\leq \bar{\gamma}$, where $(\bar{P},\bar{X},\bar{R},\bar{\alpha},\bar{\gamma})$ is the solution to Problem \eqref{ferrante minimal gamma}.

\subsection{Problem Formulation}
In this paper, we conduct a comprehensive analysis of the optimal structured and sparse $H_\infty$ controller design, and three typical problems are considered separately. To facilitate discussions, we define the stability region $\mathcal{F} \triangleq \{K\in \mathbb{R}^{n_u\times n_x}|\ A+BK\ \text{is Hurwitz} \ (\text{i.e.}, \mathrm{Re}(\lambda)<0\ \text{for all eigenvalues}\  \lambda\  \text{of}\ A+BK$)\}


\begin{myproblem}(\textit{Sparse controller design with bounded $H_\infty$ performance})\label{sparse controller design with bounded hinf performance problem}
\begin{equation}\label{p1eq}
    \mathop{\min}_{K\in \mathcal{F}} ||K||_1 \quad \text{s.t.} \quad ||T_{zd}(s)||_\infty \le \gamma,
\end{equation}
where $\gamma$ represents the performance threshold to be satisfied.
\end{myproblem}

\begin{myproblem}(\textit{Sparsity-promoting $H_\infty$ controller design})\label{sparsity-promoting hinf controller design problem}
\begin{equation}\label{p2eq}
    \mathop{\min}_{K\in \mathcal{F}} ||T_{zd}(s)||_\infty + \lambda||K||_1,
\end{equation}
where $\lambda\ge 0$ is a tuning parameter that represents the level of our emphasis on the network sparsity. 
\end{myproblem}

\begin{myrem}
    The communication costs of channels could be different
in practice. Our approach can be extended to solve this case by incorporating a weighted $l_1$ norm
which assigns a different weight to each channel based on its cost. $||K||_1$ in \eqref{p1eq} and \eqref{p2eq} can be generalized into $||W \circ K||_1$, where $\circ$ is the element-wise product and $W_{ij}$ represents the communication cost
between controller $[u(t)]_i$ and state $[x(t)]_j$. For simplicity, we assume that $W = I$ in this paper, i.e., equal
costs across all channels.
\end{myrem}
\begin{myproblem}(\textit{Structured $H_\infty$ controller design})\label{structured hinf controller design problem}
    \begin{equation}
    \mathop{\min}_{K\in \mathcal{F}} ||T_{zd}(s)||_\infty \quad \text{s.t.}\quad K\circ I_{\mathcal{S}^c} = 0.
\end{equation}
\end{myproblem}
These three problems are representative and encompass the majority of sparse controller design issues encountered in practical industrial applications. For the rest of this paper, we will analyze and solve them separately.
\section{Main Results}\label{Main Results}
\subsection{Sparse Controller Design with Bounded $H_\infty$ Performance}
In this subsection, we design a satisfactory sparse $H_\infty$ controller, ensuring that its $H_\infty$ performance remains within a specified threshold. In other words, our goal is to design the sparsest controller while simultaneously satisfying the $H_\infty$ performance constraint.
By leveraging \cref{h infinity lemma}, \cref{sparse controller design with bounded hinf performance problem} is equivalent to the following problem.

\begin{subequations}\label{equivalent DSHGT}
\begin{align}
\mathop{\min}_{K,P} &\ ||K||_1\\
\mathrm{ s.t. } &\ P \succ 0, \\
& \left[ \begin{matrix}\centering
	       &\mathrm{Sym}(P(A+BK))		&PG		&(C+DK)^\top		\\
	       &G^\top P		&-\gamma I		&H^\top		\\
	       &C+DK		&H		&-\gamma I		\\
            \end{matrix} \right] \preceq 0.
\end{align}
\end{subequations}


In the absence of the sparsity regulation term in the objective function, Problem \eqref{equivalent DSHGT} turns out to be a feasibility problem and can be solved by a standard change of variable \cite{duan2013lmis}. However, this technique is no longer viable when the sparsity regulation term is the objective we want to minimize
since this will lead to minimizing $||YX^{-1}||_1$, which is intractable. Instead of directly solving Problem \eqref{equivalent DSHGT}, we propose a linearized version of Problem \eqref{equivalent DSHGT} and clarify their relationship. We also describe how this relationship addresses the difficulties and motivates us to design a viable algorithm. At the end of this subsection, we design an ILMI algorithm to solve Problem \eqref{equivalent DSHGT} numerically. We show that this algorithm is guaranteed to converge and every limit point is a stationary point of Problem \eqref{equivalent DSHGT}.

To facilitate further discussions, we first define some critical functions as the prerequisites. Define
\begin{equation}
    f(K,P):=\frac{1}{2}(A+BK-P)^\top(A+BK-P)
\end{equation}
and its linearization around $(\Tilde{K},\Tilde{P})$ in \eqref{linearization}.
\begin{equation}\label{linearization}
\begin{aligned}
    &L_{f}(K,P;\Tilde{K},\Tilde{P}):= f(\Tilde{K},\Tilde{P})\\
    &+\frac{1}{2}\left[B(K-\Tilde{K})-(P-\Tilde{P})\right]^\top(A+B\Tilde{K}-\Tilde{P})\\
    &+\frac{1}{2}(A+B\Tilde{K}-\Tilde{P})^\top\left[B(K-\Tilde{K})-(P-\Tilde{P})\right].
\end{aligned}
\end{equation}
We further propose the following problem, which will be incorporated into our final algorithm.
\begin{subequations}\label{relaxed problem}
\begin{align}
\mathop{\min}_{K,P} &\  ||K||_1  \\
\mathrm{ s.t. } &\ P \succ 0, \label{remain constraint general}\\
& \left[
  \begin{matrix}
    &-L_f(K,P;\Tilde{K},\Tilde{P}) & * &* & * \\ 
    &\frac{1}{\sqrt{2}}(A+BK+P)  & -I &*           & *  \\ 
    &G^\top P  & 0 &-\gamma I  & *\\ 
    &C+DK & 0 &H & -\gamma I\\
  \end{matrix}\right] \preceq 0. \label{relaxed LMI}
\end{align}
\end{subequations}

It is worth noting that Problem \eqref{relaxed problem} is a convex problem and can be solved using commercial solvers. In what follows, 
we will present a result (\cref{relax theorem}) that characterizes the 
relationship between Problem \eqref{relaxed problem} and 
Problem \eqref{equivalent DSHGT}. We will then describe our 
inspiration to design the final algorithm based on 
Problem \eqref{relaxed problem}.
\begin{mythr}\label{relax theorem}
    Every feasible solution $(\bar{K},\bar{P})$ to Problem \eqref{relaxed problem} is a feasible solution to Problem \eqref{equivalent DSHGT}. For every feasible solution $(\hat{K},\hat{P})$ to Problem \eqref{equivalent DSHGT}, take $\Tilde{K}=\hat{K}, \Tilde{P}=\hat{P}$, then at least $(\hat{K},\hat{P})$ is a feasible solution to Problem \eqref{relaxed problem}.
\end{mythr}
\quad \textit{Proof}: Since all the bi-linearity of Problem \eqref{equivalent DSHGT} comes from the $(A+BK)^\top P+P(A+BK)$, it is natural to linearize this term. We apply a well-known equality and reformulate this term as follows
\begin{equation*}
\begin{split}
&(A+BK)^\top P+P(A+BK) \\
&=\frac{1}{2}(A+BK+P)^\top(A+BK+P)\\
&-\frac{1}{2}(A+BK-P)^\top(A+BK-P).
\end{split}
\end{equation*}
Consider the last term
\begin{equation*}
    \begin{split}
        &\frac{1}{2}(A+BK-P)^\top(A+BK-P)\\
        &=\frac{1}{2}\left[A+B\Tilde{K}-\Tilde{P}+B(K-\Tilde{K})-(P-\Tilde{P})\right]^\top\\
        &\left[A+B\Tilde{K}-\Tilde{P}+B(K-\Tilde{K})-(P-\Tilde{P})\right]\\
        &=\frac{1}{2}(A+B\Tilde{K}-\Tilde{P})^\top(A+B\Tilde{K}-\Tilde{P})\\
        &+\frac{1}{2}\cdot\mathrm{Sym}\left(\left[B(K-\Tilde{K})-(P-\Tilde{P})\right]^\top(A+B\Tilde{K}-\Tilde{P})\right)\\
        &+\frac{1}{2}\left[B(K-\Tilde{K})-(P-\Tilde{P})\right]^\top\left[B(K-\Tilde{K})-(P-\Tilde{P})\right]\\
        & \succeq \frac{1}{2}(A+B\Tilde{K}-\Tilde{P})^\top(A+B\Tilde{K}-\Tilde{P})\\
        &+\frac{1}{2}\cdot\mathrm{Sym}\left(\left[B(K-\Tilde{K})-(P-\Tilde{P})\right]^\top(A+B\Tilde{K}-\Tilde{P})\right)\\
        &= L_f(K,P;\Tilde{K},\Tilde{P}).
    \end{split}
\end{equation*}
The inequality follows from the semi-definiteness of the second-order incremental matrix. Thus,
\begin{equation*}
    \begin{split}
        &(A+BK)^\top P+P(A+BK)\\
        &\preceq \frac{1}{2}(A+BK+P)^\top(A+BK+P)\\
        &-L_f(K,P;\Tilde{K},\Tilde{P}).
    \end{split}
\end{equation*}
Furthermore,
\begin{equation}\label{3*3 inequality}
\begin{split}
    &\left[ \begin{matrix}\centering
	       &(A+BK)^\top P+P(A+BK)		&PG		&(C+DK)^\top		\\
	       &G^\top P		&-\gamma I		&H^\top		\\
	       &C+DK		&H		&-\gamma I		\\
            \end{matrix} \right]\\
    & \preceq \\
    &\left[
  \begin{matrix}
    &\frac{1}{2}(A+BK+P)^\top(A+BK+P)-L_f(K,P;\Tilde{K},\Tilde{P})		\\ 
    &G^\top P		\\ 
     &C+DK		\\ 
  \end{matrix}\right.                
\\
  &~~~~~~~~~~~~~~~~~~~~~~~~~~~~~~~~~~~~~~~~~~\left.
  \begin{matrix}
    &PG&(C+DK)^\top\\ 
    &-\gamma I&H^\top          \\ 
    &H&-\gamma I\\ 
  \end{matrix}\right].
\end{split}
\end{equation}
The negative definiteness of the right-hand side is equivalent to \eqref{relaxed LMI} by applying the Schur complement. Consequently, the feasible region $(K,P)$ of Problem \eqref{relaxed problem} is a subset of that of Problem \eqref{equivalent DSHGT} and we prove every feasible solution $(\bar{K},\bar{P})$ to Problem \eqref{relaxed problem} is a feasible solution to Problem \eqref{equivalent DSHGT}. Note that the inequality in \eqref{3*3 inequality} turns out to be an equality when $K=\Tilde{K}, P=\Tilde{P}$ because in this case, the linearization introduces no error. Then for every feasible solution $(\hat{K},\hat{P})$ to Problem \eqref{equivalent DSHGT}, take $\Tilde{K}=\hat{K}, \Tilde{P}=\hat{P}$, then at least $(\hat{K},\hat{P})$ is a feasible solution to Problem \eqref{relaxed problem} and possibly a convex set containing $(\hat{K},\hat{P})$ is the feasible region of Problem \eqref{relaxed problem}. $\hfill \square$

\cref{relax theorem} indicates that instead of directly solving Problem \eqref{equivalent DSHGT}, we can alternatively solve a convex counterpart (Problem \eqref{relaxed problem}), and the solution will always remain within the feasible domain of Problem \eqref{equivalent DSHGT}. Furthermore, if we consider solving Problem \eqref{relaxed problem} iteratively with linearization around the solution from the previous step, the solution is at least not worse than the previous step. This is because Problem \eqref{relaxed problem} is convex and the global minimum is not worse than the linearization point. The analysis above motivates us to design an iterative algorithm to gradually approach the stationary point, which is the highest pursuit in non-convex optimization problems.

It is worth noting that although iteratively solving Problem \eqref{relaxed problem} ensures a non-increasing sequence of objective values, this does not guarantee convergence of the algorithm. This is because the algorithm may oscillate between two points with the same objective value. To guarantee the convergence of the algorithm, the objective function must exhibit a sufficient decrease property. Therefore, we slightly modify Problem \eqref{relaxed problem} by incorporating proximal terms in the objective function, as shown in Problem \eqref{proximal relax problem}. These proximal terms avoid the oscillation issue and enforce the algorithm to gradually converge.

\begin{subequations} \label{proximal relax problem}
\begin{align}
    \mathop{\min}_{K,P}&\ ||K||_1+||K-\Tilde{K}||^2_F+||P-\Tilde{P}||^2_F\\
    \mathrm{s.t.} &\ \eqref{remain constraint general}, 
    \eqref{relaxed LMI},
\end{align}
\end{subequations}

where $\Tilde{K}, \Tilde{P}$ represent the optimization variables from the previous iteration, which will be clarified in the final algorithm. The full algorithm is shown in \cref{DSHGT algorithm}.
\begin{algorithm}[t]
\caption{Sparse $H_\infty$ controller design given $\gamma$}\label{DSHGT algorithm}

\textbf{Output:} $K^*_{A1}$;\\
Initialize  $k=0$;\\
Initialize $K_0=Y^*{X^*}^{-1}, P_0={X^*}^{-1}$ as the optimal centralized solution to Problem \eqref{unconstrained h infinity};\\
\Repeat{$||K_k-K_{k-1}||_F < \epsilon, ||P_k-P_{k-1}||_F < \epsilon$\vspace{5pt}}{
  - Solve Problem \eqref{proximal relax problem}, $\Tilde{K}=K_k, \Tilde{P}=P_k$;\\
  - Assign the solution to $K_{k+1}, P_{k+1}$;\\
  - $k = k + 1$;}
\Return{$K_k$};
\end{algorithm}

Now we start to evaluate the optimality of \cref{DSHGT algorithm}. To facilitate discussions on the optimality, we first provide the Lagrangian of Problem \eqref{equivalent DSHGT} and then give the Karush-Kuhn-Tucker (KKT) conditions of Problem \eqref{equivalent DSHGT} and the definition of the stationary point.

The Lagrangian of Problem \eqref{equivalent DSHGT} is shown as
\begin{equation*}
    L(K,P,\Lambda_1,\Lambda_2) = ||K||_1 +\langle\Lambda_1,N(K,P) \rangle - \langle\Lambda_2,P\rangle,
\end{equation*}
where $\Lambda_1,\Lambda_2$ are the Lagrangian multipliers, and we denote 
\begin{equation}
\label{N}
    N(K,P):= \left[ \begin{matrix}\centering
	       &\mathrm{Sym}(P(A+BK))		&PG		&(C+DK)^\top		\\
	       &G^\top P		&-\gamma I		&H^\top		\\
	       &C+DK		&H		&-\gamma I		\\
            \end{matrix} \right]
\end{equation}
for simplicity. Therefore, we can derive the KKT conditions as follows
\begin{subequations}\label{KKT of equivalent DSHGT}
\begin{align}
        &\Lambda^*_1 \succeq 0, \Lambda^*_2 \succeq 0,\\
        &N(K^*,P^*)\preceq 0, P^*\succeq 0,\\
        &\langle \Lambda^*_1, N(K^*,P^*) \rangle = 0, \langle \Lambda^*_2,P \rangle = 0,\\
        &\partial ||K^*||_1 + \nabla_K \langle \Lambda^*_1, N(K^*,P^*) \rangle = 0,\\
        &\nabla_P\langle \Lambda^*_1, N(K^*,P^*) \rangle - \Lambda_2 = 0,&
\end{align}
\end{subequations}
where $\Lambda^*_1, \Lambda^*_2$ are the KKT multipliers. With the KKT conditions, we define the stationary point below.

\begin{mydef}
    $(K,P)$ is called a stationary point of Problem \eqref{equivalent DSHGT} if it satisfies the corresponding KKT conditions \eqref{KKT of equivalent DSHGT}.
\end{mydef}

\begin{myrem}
    Generally speaking, a stationary point is not necessarily a locally optimal point because of the existence of saddle points. Although it is not a sufficient condition for local optimality, we can still remove points that are not locally optimal if we can show they do not satisfy the KKT conditions. The stationary point is in general the highest pursuit in non-convex optimization problems.
\end{myrem}

\begin{mythr}\label{convergence theorem}
    \cref{DSHGT algorithm} generates a solution sequence that has at least one limit point, and each limit point is a stationary point of \cref{sparse controller design with bounded hinf performance problem}.
\end{mythr}
\quad \textit{Proof}: Since Problem \eqref{proximal relax problem} is a convex problem, we can obtain the global optimal solution through commercial solvers. From the structure of \cref{DSHGT algorithm}, we immediately obtain the following relationship at each step
    \begin{equation}\label{sufficient descrease}
        ||K_{k+1}||_1 + ||K_{k+1}-K_k||^2_F+||P_{k+1}-P_k||^2_F
        \le  ||K_k||_1
    \end{equation}
    because the global optimal point is at least not worse than the linearization point. From \eqref{sufficient descrease}, we further obtain if $(K_{k+1},P_{k+1})$ is not equal to $ (K_k,P_k)$, then the objective function will decrease. Since the objective function is lower-bounded, the objective value must finally converge. From the relation
    \begin{equation*} 
    \begin{split}
    &\mathop{\lim}_{k\rightarrow \infty} ||K_{k+1}||_1 + ||K_{k+1}-K_k||^2_F+||P_{k+1}-P_k||^2_F\\
    & \le \mathop{\lim}_{k\rightarrow \infty} ||K_k||_1,
    \end{split}
    \end{equation*}
    we have
    \begin{equation*}
        \mathop{\lim}_{k\rightarrow \infty}||K_{k+1}-K_k||^2_F+||P_{k+1}-P_k||^2_F=0.
    \end{equation*}
    Therefore the solution sequence must have at least one limit point. According to \cref{DSHGT algorithm}, the solution sequence is obtained by iteratively solving Problem \eqref{proximal relax problem}. Thus, it is natural to consider the optimality of Problem \eqref{proximal relax problem} at each iteration. After recovering \eqref{relaxed LMI} using the Schur complement, we provide the Lagrangian of Problem \eqref{proximal relax problem} at the $k+1$ iteration.
    \begin{equation*}
    \begin{split}
        &L(K,P,[\Gamma_1]_{k+1},[\Gamma_2]_{k+1}) \\
        &=  ||K||_1 +||K-K_k||^2_F+||P-P_k||^2_F\\
        &+ \langle [\Gamma_1]_{k+1},M(K,P;K_k,P_k) \rangle - \langle [\Gamma_2]_{k+1},P \rangle,\\
    \end{split}
    \end{equation*}
    where $[\Gamma_1]_{k+1},[\Gamma_2]_{k+1}$ denotes the Lagrangian multipliers, and we denote 
\begin{equation}
\label{M}
\begin{split}
    &M(K,P;K_k,P_k) = \\
&\left[
  \begin{matrix}
    &\frac{1}{2}(A+BK+P)^\top(A+BK+P)-L_f(K,P;K_k,P_k)		\\ 
    &G^\top P		\\ 
     &C+DK		\\ 
  \end{matrix}\right.                
\\
  &~~~~~~~~~~~~~~~~~~~~~~~~~~~~~~~~~~~~~~~~~~~~\left.
  \begin{matrix}
    &PG&(C+DK)^\top\\ 
    &-\gamma I&H^\top          \\ 
    &H&-\gamma I\\ 
  \end{matrix}\right]
\end{split}
\end{equation}
for conciseness. Since Problem \eqref{proximal relax problem} is a convex problem, the optimal solution $(K_{k+1},P_{k+1})$ must satisfy the corresponding KKT conditions. The KKT conditions of Problem \eqref{proximal relax problem} at the $k+1$ iteration are listed in \eqref{KKT of proximal relax problem}, where
$[\Gamma^*_1]_{k+1},[\Gamma^*_2]_{k+1}$ are known  as the KKT multipliers. Besides, the equality $N(K,P)=M(K,P;K,P)$ holds according to their definitions \eqref{N}, \eqref{M}. Then we take the limit for all KKT conditions \eqref{KKT of proximal relax problem} and use subscript $\infty$ to denote any limit point. Let $\Lambda^*_1=[\Gamma^*_1]_\infty$, $\Lambda^*_2=[\Gamma^*_2]_\infty$. The limit point $(K_\infty,P_\infty)$ satisfies the KKT conditions of \eqref{KKT of equivalent DSHGT}, which means that $(K_\infty,P_\infty)$ is a stationary point of Problem \eqref{equivalent DSHGT} and thus \cref{sparse controller design with bounded hinf performance problem}.$\hfill \square$

\begin{figure*}[t]
    \begin{subequations}\label{KKT of proximal relax problem}
    \begin{align}
        &[\Gamma^*_1]_{k+1} \succeq0, [\Gamma^*_2]_{k+1} \succeq0,\\
        &M(K_{k+1},P_{k+1};K_k,P_k)\preceq0,  P_{k+1}\succeq0, \\
        &\langle [\Gamma^*_1]_{k+1}, M(K_{k+1},P_{k+1};K_k,P_k) \rangle = 0, \langle [\Gamma^*_2]_{k+1},P_{k+1} \rangle =0, \\
         &\partial||K_{k+1}||_1 + (K_{k+1}-K_k)+\nabla_K\langle [\Gamma^*_1]_{k+1}, M(K_{k+1},P_{k+1};K_k,P_k)\rangle=0,\\
         &(P_{k+1}-P_k)+\nabla_P \langle [\Gamma^*_1]_{k+1}, M(K_{k+1},P_{k+1};K_k,P_k) \rangle - [\Gamma^*_2]_{k+1}=0,&
    \end{align}
\end{subequations}
\end{figure*}
This subsection provides a reliable algorithm with guaranteed convergence to solve \cref{sparse controller design with bounded hinf performance problem}, a widely recognized scenario in sparse controller design. The next subsection will focus on the trade-off between system performance and communication burdens.
\subsection{Sparsity-promoting $H_\infty$ Controller Design}
Different from the assumption of the previous subsection, in reality, sometimes the performance threshold is unknown and we are interested in the relationship between system performance and topology sparsity, which motivates us to consider \cref{sparsity-promoting hinf controller design problem}. By utilizing \cref{h infinity lemma}, \cref{sparsity-promoting hinf controller design problem} is equivalent to the following problem.

\begin{subequations}\label{equivalent SPHCD}
\begin{align}
            \mathop{\min}_{K,P,\gamma} &\ \gamma+ \lambda||K||_1 \\
            \mathrm{ s.t. } &\ P \succ 0 ,\\
            & \left[ \begin{matrix}\centering
	       &\mathrm{Sym}(P(A+BK))		&PG		&(C+DK)^\top		\\
	       &G^\top P		&-\gamma I		&H^\top		\\
	       &C+DK		&H		&-\gamma I		\\
            \end{matrix} \right] \preceq 0 .
\end{align}
\end{subequations}

\begin{algorithm}[t]
\caption{Sparsity promoting $H_\infty$ controller design}\label{SPHCD algorithm}
\textbf{Output:} $K^*_{A2}, P^*_{A2}, \gamma^*_{A2}$;\\
Initialize  $k=0$;\\
Initialize $K_0=Y^*{X^*}^{-1}, P_0={X^*}^{-1},\gamma_0=\gamma^*$ as the optimal centralized solution to Problem \eqref{unconstrained h infinity};\\
\Repeat{$||K_k-K_{k-1}||_F < \epsilon, ||P_k-P_{k-1}||_F < \epsilon$}{
  - Solve Problem \eqref{proximal SPHCD} with $\Tilde{K}=K_k,\Tilde{P}=P_k$ ;\\
  - Assign the solutions to $K_{k+1}, P_{k+1}, \gamma_{k+1}$;\\
  - $k = k + 1$;
}
\Return{$K_k,P_k,\gamma_k$};
\end{algorithm}

Different from Problem \eqref{equivalent DSHGT}, here we incorporate $\gamma$ into the objective function as a variable to optimize. The parameter $\lambda$ is the penalty weight that represents the level of our emphasis on sparsity. As $\lambda$ increases, the solution $K$ would become sparser. 

Note that the LMI constraint of Problem \eqref{equivalent SPHCD} also contains the bi-linear term. To tackle this difficulty, we use the same technique as the previous section: i) applying linearization; (ii) leveraging the Schur complement to transform the bi-linear matrix inequality into an LMI. The problem that will be adopted is shown in Problem \eqref{proximal SPHCD}, where we also incorporate the proximal terms to guarantee sufficient decreasing. The $\Tilde{K},\Tilde{P}$ are the linearization point. The full algorithm is shown in \cref{SPHCD algorithm}.

\begin{mythr}
    \cref{SPHCD algorithm} generates a solution sequence that has at least one limit point, and each limit point is a stationary point of \cref{sparsity-promoting hinf controller design problem}
\end{mythr}
\quad \textit{Proof}: It is similar to \cref{convergence theorem} and is omitted.$\hfill \square$

This subsection provides a similar algorithm to find the controller structure with a convergence guarantee. However, the solution is not optimal and can be improved due to the $l_1$ relaxation. Because the $l_1$ regulation penalizes elements with large magnitude, the solution returned tends to unnecessarily have small elements within the nonzero pattern. This deviates from our ultimate goal, which is penalizing the matrix cardinality instead of value. Therefore, we propose to refine our solution by solving a structured controller design problem over the structure obtained in this part. 
\begin{subequations}    \label{proximal SPHCD}
\begin{align}
    \mathop{\min}_{K,P,\gamma} &\ \gamma+\lambda ||K||_1+||K-\Tilde{K}||^2_F+||P-\Tilde{P}||^2_F\\
    \mathrm{ s.t. } &\ P \succ 0 \label{SPHCD general},\\
& \left[
  \begin{matrix}
    &-L_f(K,P;\Tilde{K},\Tilde{P}) & * &* & * \\ 
    &\frac{1}{\sqrt{2}}(A+BK+P)  & -I &*           & *  \\ 
    &G^\top P  & 0 &-\gamma I  & *\\ 
    &C+DK & 0 &H & -\gamma I \\
  \end{matrix}\right] \preceq 0 \label{SPHCD relaxed LMI}.
\end{align}
\end{subequations}

\begin{algorithm}[t]
\caption{Structured $H_\infty$ controller design}\label{structured algorithm}
\textbf{Output:} $K^*_{A3}, P^*_{A3}, \gamma^*_{A3}$;\\
Initialize  $k=0$;\\
Solve Probelm \eqref{ferrante minimal gamma} and denote the solution as $(\bar{P},\bar{X},\bar{R},\bar{\alpha},\bar{\gamma})$;\\
Initialize $K_0=\bar{R}{\bar{X}}^{-1}, P_0=\bar{P}^{-1}, \gamma_0=\bar{\gamma}$ ;\\
\Repeat{$||K_k-K_{k-1}||_F < \epsilon, ||P_k-P_{k-1}||_F < \epsilon$\vspace{5pt}}{
  - Solve Problem \eqref{structured problem} with $\Tilde{K}=K_k, \Tilde{P}=P_k$;\\
  - Assign the solution to $K_{k+1}, P_{k+1}, \gamma_{k+1}$;\\
  - $k = k + 1$;
}
\Return{$K_k,P_k,\gamma_k$};
\end{algorithm}

\subsection{Structured $H_\infty$ Controller Design}
Different from the previous two subsections, we now consider designing an $H_\infty$ controller with a specified structural pattern $K\in \mathcal{S}$, where $\mathcal{S}$ is a linear subspace. This is motivated by the reality that a controller's structure is sometimes naturally or artificially restricted, and a controller lacking such structure may fail in practice. For example, the controller for a micro-electro-mechanical-system (MEMS) must adhere to a specific structural pattern due to the internal electrical circuitry \cite{fan2008structured}. 

Ferrante et al. \cite{ferrante2020lmi} proposed a viable approach to address this problem, as described in \cref{ferrante results}. However, there are some limitations: (i) their approach explores a structured controller within a convex subset of the nonconvex feasibility region, omitting some potential solutions and can lead to suboptimal results; (ii) the minimal value $\bar{\gamma}$ returned by their approach (Problem \eqref{ferrante minimal gamma}) is not the $H_\infty$ norm itself, but rather an upper bound of the $H_\infty$ norm correspnding to the structured controller $\bar{K}$. To clarify this point, we consider the following problem
\begin{subequations}\label{auxiliary problem}
\begin{align}
&\mathop{\min}_{P,X,\alpha,\gamma} \ \gamma \quad\mathrm{ s.t. }\\
&\  P\succ 0, \alpha>0,\\
&\ \scalebox{0.87}{$\mathrm{Sym}\left(\begin{bmatrix}
                &(A+B\bar{K})X &\alpha(A+B\bar{K})X+P &0 &G\\
                &-X  &-\alpha X  &0  &0\\
                &(C+D\bar{K})X  &\alpha(C+D\bar{K})X  &-\frac{\gamma}{2}I  &H\\
                &0  &0  &0  &-\frac{\gamma}{2}I
            \end{bmatrix}\right)\prec 0$}
\end{align}
\end{subequations}
Recall that $(\bar{P},\bar{X},\bar{R},\bar{\alpha},\bar{\gamma})$ is the solution to Problem \eqref{ferrante minimal gamma} and $\bar{K}=\bar{R}\bar{X}^{-1}$. Consequently, it is straightforward to see that $(\bar{P},\bar{X},\bar{\alpha},\bar{\gamma})$ is a feasible solution to Problem \eqref{auxiliary problem}. Thus $\bar{\gamma}$ serves as an upper bound of $\gamma^*$, which is returned by Problem \eqref{auxiliary problem}. According to \cref{equivalence}, $\bar{\gamma}$ also upper bounds the $H_\infty$ norm corresponding to $\bar{K}$.

On the contrary, our approaches can effectively address these limitations accordingly: (i) our ILMI algorithm can iteratively consider different convex subsets within the nonconvex feasibility domain and refine the solution at each step. This allows our method to explore a larger feasibility domain than that proposed in \cite{ferrante2020lmi}; (ii) denoting the output of our algorithm as $(K^*_{A3}, P^*_{A3}, \gamma^*_{A3})$, $\gamma^*_{A3}$ is exactly the $H_\infty$ norm corresponding to $K^*_{A3}$. 

Despite the limitations mentioned earlier, we can leverage the solution returned by Problem \eqref{ferrante minimal gamma} as the initial point of our ILMI algorithm discussed later for faster convergence, since this solution is closer to the optimal solution than a randomly chosen feasible one.

Following the concepts from the previous two subsections, we introduce Problem \eqref{structured original problem} and Problem \eqref{structured problem}. Our goal is to solve Problem \eqref{structured original problem}, which is a structured optimal controller design problem. However, since Problem \eqref{structured original problem} is non-convex, we iteratively solve Problem \eqref{structured problem} instead. Finally, we demonstrate that the solution returned by our algorithm is a stationary point of Problem \eqref{structured original problem}. 
\begin{subequations}\label{structured original problem}
\begin{align}
            \mathop{\min}_{K,P,\gamma} &\ \gamma \\
             \mathrm{ s.t. } &\ P \succ 0 , K \circ I_{\mathcal{S}^c} = 0 ,\\
            & \left[ \begin{matrix}\centering
	       &\mathrm{Sym}(P(A+BK))		&PG		&(C+DK)^\top		\\
	       &G^\top P		&-\gamma I		&H^\top		\\
	       &C+DK		&H		&-\gamma I		\\
            \end{matrix} \right] \preceq 0 .
\end{align}
\end{subequations}

\begin{subequations}\label{structured problem}
\begin{align}
        \mathop{\min}_{K,P,\gamma}&\ \gamma+||K-\Tilde{K}||^2_F+||P-\Tilde{P}||^2_F\\
        \text{ s.t. } &\ P \succ 0 , K \circ I_{\mathcal{S}^c} = 0 ,\\
         & \left[
  \begin{matrix}
    &-L_f(K,P;\Tilde{K},\Tilde{P}) & * &* & * \\ 
    &\frac{1}{\sqrt{2}}(A+BK+P)  & -I &*           & *  \\ 
    &G^\top P  & 0 &-\gamma I  & *\\ 
    &C+DK & 0 &H & -\gamma I \\
  \end{matrix}\right] \preceq 0 .
\end{align}
\end{subequations}

The full algorithm is shown in \cref{structured algorithm}.
\begin{mythr}
    \cref{structured algorithm} generates a solution sequence that has at least one limit point, and each limit point is a stationary point of Problem \eqref{structured original problem}.
\end{mythr}
\quad \textit{Proof}: It is similar to \cref{convergence theorem} and is omitted.
$\hfill \square$

\begin{myrem}
    Apart from increasing convergence speed, another reason for utilizing the solution returned by Problem \eqref{ferrante minimal gamma} for initialization is to avoid infeasibility in the first iteration. Given $(\bar{P},\bar{X},\bar{R},\bar{\alpha},\bar{\gamma})$, we compute $\Tilde{K}=\bar{R}{\bar{X}}^{-1}, \Tilde{P}=\bar{P}^{-1}, \Tilde{\gamma}=\bar{\gamma}$, where $(\Tilde{K}, \Tilde{P}, \Tilde{\gamma})$ are variables specified in Problem \eqref{structured problem}. At least $(\Tilde{K},\Tilde{P},\Tilde{\gamma})$, and potentially other superior solutions, are feasible for Problem \eqref{structured problem} since $\Tilde{K}\circ I_{\mathcal{S}^c}=0$ holds directly and $\Tilde{\gamma}$ upper bounds the $H_\infty$ norm corresponding to $\Tilde{K}$. However, this result does not necessarily hold for a general $(\Tilde{K},\Tilde{P},\Tilde{\gamma})$ where $\Tilde{K}\circ I_{\mathcal{S}^c}\neq0$. This is because in Problem \eqref{structured problem}, we consider a convex subset of $(K,P)$ around $(\Tilde{K},\Tilde{P})$ within the nonconvex feasibility region. We cannot guarantee that there exists a $K$ in this region that satisfies $\Tilde{K}\circ I_{\mathcal{S}^c}=0$.
\end{myrem}

\section{Simulations}\label{Simulations}
In this section, we provide several numerical simulations to verify our results. The numerical example for \cref{sparse controller design with bounded hinf performance problem} and \cref{sparsity-promoting hinf controller design problem} is described as follows: consider the mass-spring system with $N$ masses on a line \cite{lin2011augmented}. The dynamic system can be represented as $x_1=[p_1,\ldots,p_N]^\top$ and $x_2 = \dot{x_1}$, where $p_i$ is the displacement of the $i$-th mass from its reference point. The state-space model can be modeled as \eqref{original system} with 
\begin{align*}
        &A = \begin{bmatrix}
	~0&		~I~~ \\
	~T&		~0~~
\end{bmatrix}, ~B=\begin{bmatrix}
	~0~~\\
	~I~~\\
\end{bmatrix}, ~G=\begin{bmatrix}
	~0~~\\
	~I~~\\
\end{bmatrix},\\
&C = \begin{bmatrix}
	~I~~\\
	~0~~\\
\end{bmatrix},~D=\begin{bmatrix}
	~0~~\\
	~2I~~\\
\end{bmatrix},~H=\begin{bmatrix}
	~0~~\\
	~2I~~\\
\end{bmatrix}.
\end{align*}
where $T\in \mathbb{R}^{N\times N}$ is a Toeplitz matrix with -2 on its main diagonal, 1 on its first sub- and super-diagonal, and 0 elsewhere. We assume $N=20, n_x=40,n_u=20,n_d=20,n_z=60$. The numerical examples for \cref{structured hinf controller design problem} are the same as those in \cite{ferrante2020lmi}, which will be clarified later. Throughout this section, we choose $\epsilon=1e-3$. For the rest of this section, \cref{sparse controller design with bounded hinf performance problem}, \cref{sparsity-promoting hinf controller design problem}, and \cref{structured hinf controller design problem} will be considered separately.

\subsection{Sparse Controller Design with Bounded $H_\infty$ Performance}
\begin{figure}[t]
    \centering
    \hfill 
    \begin{subfigure}[b]{0.49\columnwidth}
        \includegraphics[width=\linewidth]{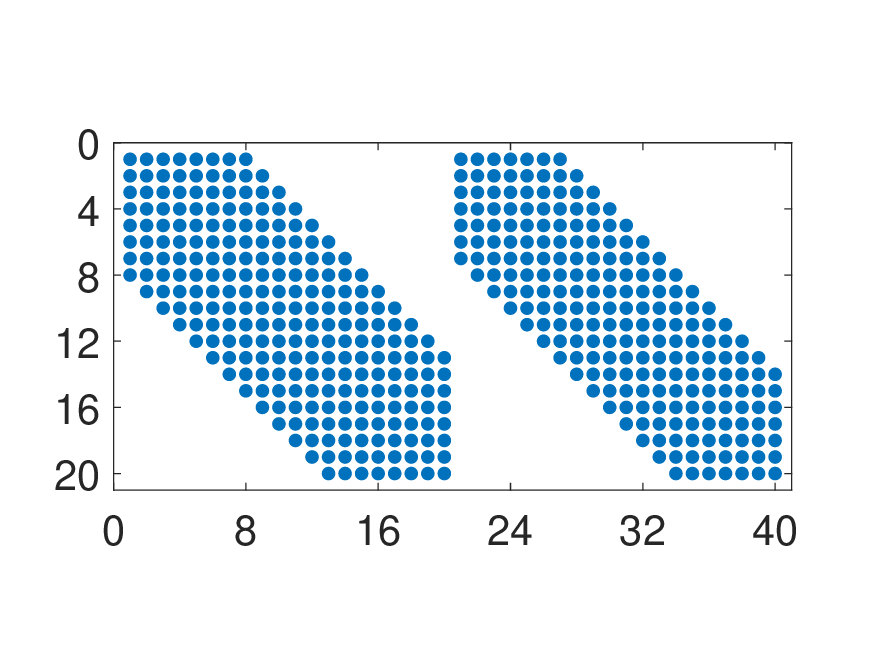}
    \end{subfigure}
    \hfill
    \begin{subfigure}[b]{0.49\columnwidth}
        \includegraphics[width=\linewidth]{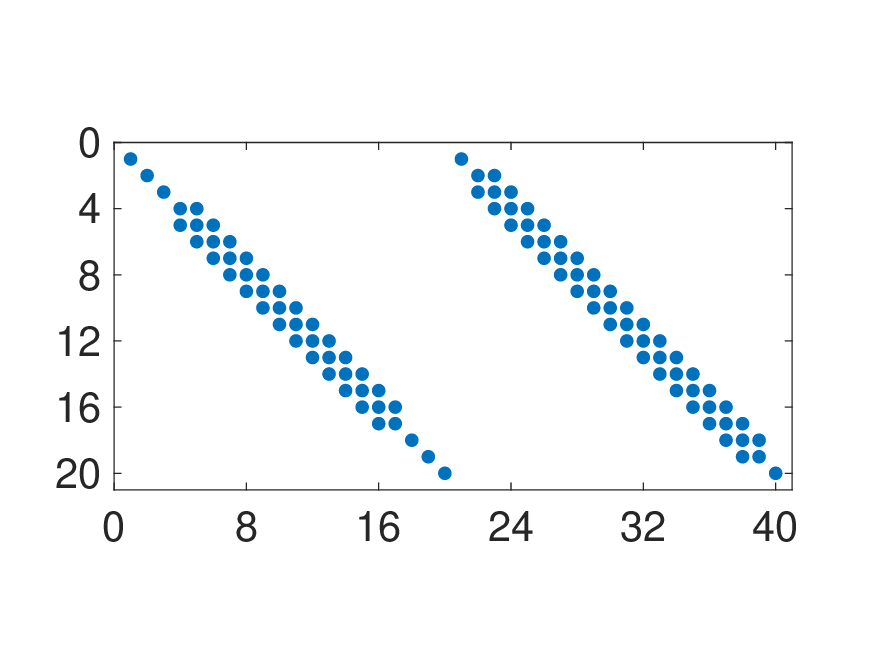}
    \end{subfigure}
    \hfill
    \begin{subfigure}[b]{0.49\columnwidth}
        \includegraphics[width=\linewidth]{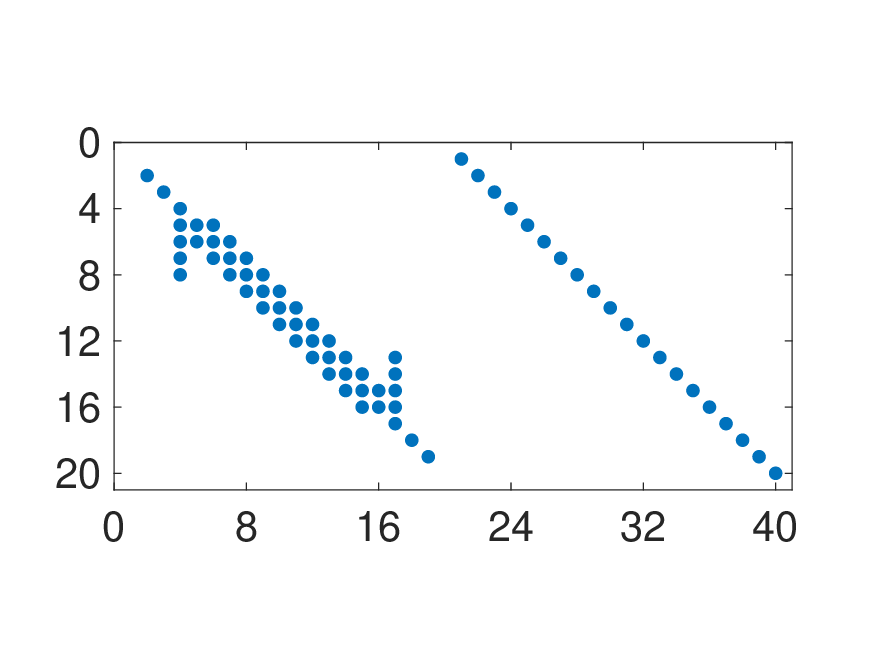}
    \end{subfigure}
    \hfill
    \begin{subfigure}[b]{0.49\columnwidth}
        \includegraphics[width=\linewidth]{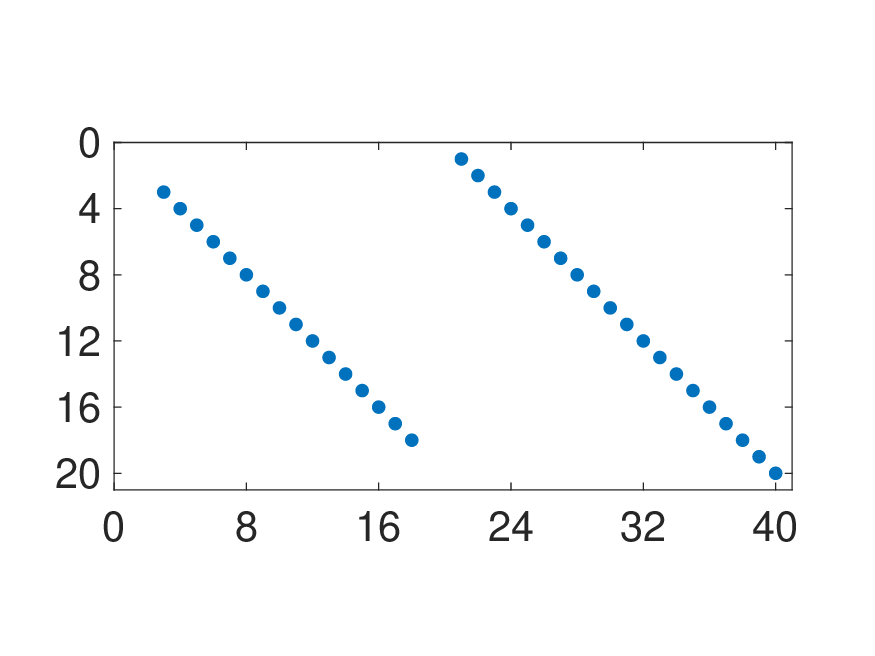}
    \end{subfigure}
    \caption{The sparsity patterns of $K_0$,$K_1$,$K_2$,$K^*_{A1}$ (from left to right, top to bottom). The nonzero elements are labeled using blue dots.}
\label{sparsity pattern}
\end{figure}

\begin{figure}[t]
        \centering   \includegraphics[width=0.44\textwidth]{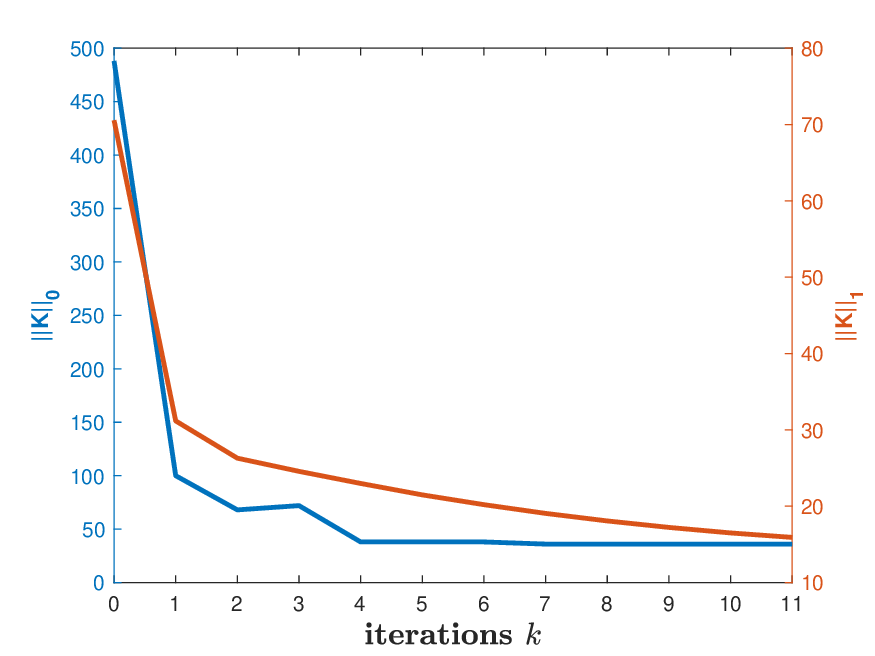}
        \caption{The evolution of $||K||_0$ and $||K||_1$.}
        \label{evolution ||K||_0 ||K||_1}
\end{figure}

By solving Problem \eqref{unconstrained h infinity}, we obtain the global optimal $H_\infty$ controller $K^*=Y^*{X^*}^{-1}$ and global minimal $H_\infty$ norm $\gamma^*=2$. We require the $H_\infty$ norm should not be larger than $\gamma=5$. We run \cref{DSHGT algorithm} to compute the solution. The solution patterns at several starting iterations are shown in \cref{sparsity pattern}. It is shown that the number of nonzero elements declines rapidly. The algorithm takes $4$ iterations to reach the optimal solution with 38 nonzero elements in 28.5563 seconds. The evolution of $||K||_0$ and $||K||_1$ are shown in \cref{evolution ||K||_0 ||K||_1}. The monotonic decreasing $||K||_1$ matches our theoretical results. However, since $||K||_1$ is an approximate of $||K||_0$, $||K||_0$ is not guaranteed to monotonically decrease. 

\begin{figure}[t]
        \centering

        \includegraphics[width=0.44\textwidth]{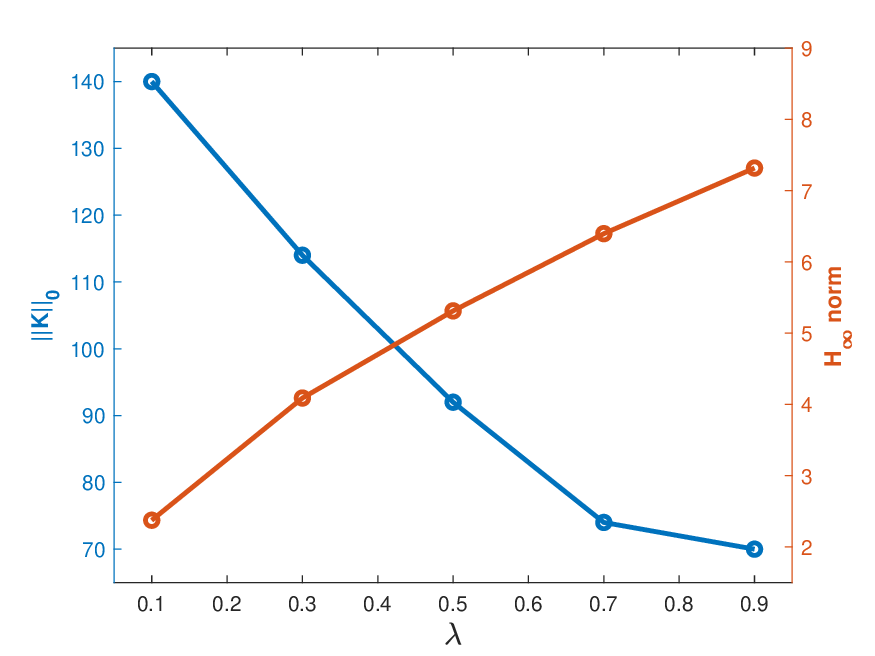}
        \caption{The relationship between $||K||_0$, $H_\infty$ norm and $\lambda$.}
        \label{relationship versus lambda}
\end{figure}
\subsection{Sparsity-promoting $H_\infty$ Controller Design}
In this part, we illustrate the trade-off between system performance and the controller sparsity in \cref{relationship versus lambda}. For each $\lambda$, \cref{SPHCD algorithm} can return a satisfactory solution within 10 iterations and 100 seconds. When $\lambda$ is small, we impose a small penalty on the controller density, resulting in denser controllers and higher performance. As $\lambda$ increases, the penalty on density also increases, leading to sparser controllers and a corresponding degradation in system performance.

\subsection{Structured $H_\infty$ Controller Design}
In this section, we will validate the effectiveness of our methods on two numerical examples discussed in \cite{ferrante2020lmi}, compare our solutions with those obtained in \cite{ferrante2020lmi} to demonstrate our superiority. 
\subsubsection{Network Decentralized Control}
We consider a water distribution system consisting of $5$ subsystems \cite{blanchini2014network}. The overall dynamics can be written in the form of \eqref{original system} with 
\begin{equation}
    \begin{split}
        A &= \mathrm{diag}(A_1,A_2,\ldots,A_5), C = I, D = 0, H = 0,\\
        B &= \begin{bmatrix}
            B_u &-B_d &0 &0 &0 &0\\
            0 &B_u &-B_d &0 &0 &-B_d\\
            0 &0 &B_d &-B_u &0 &0\\
            0 &0 &0 &B_d &B_u &0\\
            0 &0 &0 &0 &B_u &B_d
        \end{bmatrix}, G = I.
    \end{split}
\end{equation}
where for all $i=1,2,\ldots,5$,
\begin{equation}
    A_i=\begin{bmatrix}
        -\xi_i&\beta_i&0\\
        \xi_i &-\beta_i &0\\
        0 &1 &0
    \end{bmatrix},  B_d = \begin{bmatrix}
        0\\1\\0
    \end{bmatrix}
     , B_u=\begin{bmatrix}
         1\\0\\0
     \end{bmatrix}.
\end{equation}
In particular, we select $\xi_1=15, \xi_2=20, \xi_3=16, \xi_4=16.7, \xi_5=14, \beta_1=0, \beta_2=0, \beta_3=12, \beta_4=0, \beta_5=22.$ According to \cite{blanchini2014network}, a controller is \textit{decentralized in the sense of networks} if $K$ has the same structural pattern as $B^\top$. Therefore, the structural constraint can be expressed as $K\in \mathcal{S}$, with
\begin{equation}
    I_{\mathcal{S}^c}=\scalebox{0.9}{$\left[\begin{array}{ccccccccccccccc}
        1 &1 &1 &0 &0 &0 &0 &0 &0 &0 &0 &0 &0 &0 &0\\
        1 &1 &1 &1 &1 &1 &0 &0 &0 &0 &0 &0 &0 &0 &0\\
        0 &0 &0 &1 &1 &1 &1 &1 &1 &0 &0 &0 &0 &0 &0\\
        0 &0 &0 &0 &0 &0 &1 &1 &1 &1 &1 &1 &0 &0 &0\\
        0 &0 &0 &0 &0 &0 &0 &0 &0 &1 &1 &1 &1 &1 &1\\
        0 &0 &0 &1 &1 &1 &0 &0 &0 &0 &0 &0 &1 &1 &1
    \end{array}\right]$}.
\end{equation}
We run \cref{structured algorithm} to compute a structured controller. In the initialization part, we first solve Problem \eqref{ferrante minimal gamma} and denote the solution as $(\bar{P},\bar{X},\bar{R},\bar{\alpha},\bar{\gamma})$, where $\bar{\alpha}=0.1554$ and $\bar{\gamma}=1.7887$ \cite{ferrante2020lmi}. \cref{structured algorithm} returns $K^*_{A3}, P^*_{A3}, \gamma^*_{A3}$ with $\gamma^*_{A3}=1.65<1.7887$. The results demonstrate the superiority of our methods over those proposed in \cite{ferrante2020lmi}.

\subsubsection{Overlapping Control}
We consider the control with overlapping information structure constraints \cite[Example 2.16]{zecevic2010control}. The system dynamics can be expressed in the form of \eqref{original system} with
\begin{equation}
\begin{split}
    A& = \begin{bmatrix}
        1 &4 &0\\
        1 &2 &2
    \end{bmatrix}, B = \begin{bmatrix}
        1 &0\\
        0 &0\\
        0 &1
    \end{bmatrix}, G = \begin{bmatrix}
        0\\0\\1
    \end{bmatrix},\\
    C &= \begin{bmatrix}
        0\\1\\0
    \end{bmatrix}, D = 0, H=0.
\end{split}
\end{equation}
We specify a structural constraint with
\begin{equation}
    I_{\mathcal{S}^c} = \begin{bmatrix}
        1 &1 &0\\
        0 &1 &1
    \end{bmatrix}.
\end{equation}
Same as the previous example, we run \cref{structured algorithm}. We initialize the algorithm by solving Problem \eqref{ferrante minimal gamma} and obtain the solution $(\bar{P},\bar{X},\bar{R},\bar{\alpha},\bar{\gamma})$, where $\bar{\alpha}=0.09$ and $\bar{\gamma}=0.13724$ \cite{ferrante2020lmi}. \cref{structured algorithm} returns $K^*_{A3}, P^*_{A3}, \gamma^*_{A3}$ with $\gamma^*_{A3}=0.06<0.13724$, which validate our advantages over \cite{ferrante2020lmi}.

\section{Conclusions}\label{Conclusions}
In this paper, we conducted a comprehensive analysis of the optimal sparse $H_\infty$ controller design, and three typical problems were considered separately. For each problem, we applied a novel linearization to relax the bilinear matrix inequality into an LMI, and we showed that any feasible solution to the relaxed problem was feasible for the original problem, which motivated us to develop an ILMI algorithm to compute the solution. We further characterized the first-order optimality of the original problem using KKT conditions. Moreover, by incorporating proximal terms into the objective function, we showed that our algorithm is guaranteed to converge and each limit point is a stationary point of the original problem, which is the highest pursuit in non-convex optimization problems. Finally, the effectiveness of our algorithm was validated via numerical simulations. For the non-structured design, our algorithms effectively reduced the number of non-zero elements, thereby lowering communication costs. In the case of structured design, our approaches showed significant improvement over those proposed in \cite{ferrante2020lmi}.

\footnotesize{ 
\bibliographystyle{IEEEtran}
\bibliography{reference}
}
\end{document}